\documentclass[titlepage,11pt]{article}
\usepackage{color}

\usepackage{latexsym}
\usepackage{amsfonts}
\usepackage{amsmath}
\usepackage[T1]{fontenc}
\newcommand\blackslug{\hbox{\hskip 1pt \vrule width 4pt height 8pt depth 1.5pt
        \hskip 1pt}}
\newcommand\bbox{\hfill \quad \blackslug \bigbreak}
\def\dd{\hbox{-}}

\def\ll{,\ldots,}
\def\cupcup{\cup\cdots\cup}

\title{Polynomial bounds for chromatic number.\\ IV. A near-polynomial bound for excluding the five-vertex path}
\author{Alex Scott\thanks{Research supported by EPSRC grant EP/V007327/1.}\\
Mathematical Institute, University of Oxford, Oxford OX2 6GG, UK
\\
\\
Paul Seymour\thanks{Supported by AFOSR grant
A9550-19-1-0187, and by NSF grant  DMS-1800053.}\\
Princeton University, Princeton, NJ 08544
\\
\\
Sophie Spirkl\thanks{We acknowledge the support of the Natural Sciences and Engineering Research
Council of Canada (NSERC) [funding reference number RGPIN-2020-03912].
Cette recherche a \'et\'e financ\'ee par le Conseil de recherches en sciences
naturelles et en g\'enie du Canada (CRSNG) [num\'ero de r\'ef\'erence
RGPIN-2020-03912].}\\
University of Waterloo, Waterloo, Ontario N2L3G1, Canada}

\date{August 8, 2021; revised October 2, 2022}

\newtheorem{thm}{}[section]

\newcommand{\Proof}{\noindent{\bf Proof.}\ \ }

\begin{document}
\maketitle
\begin{abstract}
	\noindent A graph $G$ is {\em $H$-free} if it 
	has no induced subgraph isomorphic to $H$.
	We prove that a $P_5$-free graph with clique number $\omega\ge 3$ has chromatic number at most  $\omega^{\log_2(\omega)}$.  
	The best previous result
	was an exponential upper bound $(5/27)3^{\omega}$,
	due to Esperet, Lemoine, Maffray, and Morel.
	A polynomial bound would imply that  the
	celebrated Erd\H{o}s-Hajnal conjecture holds for $P_5$, which is the smallest open case. Thus there is great interest in whether there is a polynomial bound for $P_5$-free graphs, and our result 
	is an attempt to approach that.
	
%	There is a general conjecture, due to Esperet, that for every forest $H$, all $H$-free graphs have chromatic number at most
%a polynomial function of their clique number. This would imply the
%	celebrated Erd\H{o}s-Hajnal conjecture for $H$ and in particular for $P_5$, which is the smallest open case of the Erd\H{o}s-Hajnal conjecture. Thus there is great interest in Esperet's conjecture for $P_5$-free graphs, and our result 
%	is an attempt to approach that.
\end{abstract}

\section{Introduction}
If $G,H$ are graphs, we say $G$ is {\em $H$-free} if no induced subgraph of $G$ is isomorphic
to $H$; and for a graph $G$, we denote the number of vertices, the chromatic number, the size of the largest clique, and the size of the largest stable set 
by $|G|, \chi(G), \omega(G),\alpha(G)$ respectively.

The $k$-vertex path is denoted by $P_k$, and $P_4$-free graphs are well-understood; every $P_4$-free graph $G$ with more than one vertex is either disconnected or disconnected
in the complement, which implies that $\chi(G)=\omega(G)$. Here we study how $\chi(G)$ depends on $\omega(G)$ for $P_5$-free graphs $G$.

The Gy\'arf\'as-Sumner conjecture~\cite{gyarfas,sumner} says:
\begin{thm}\label{GYconj}
	{\bf Conjecture: }For every forest $H$ there is a function $f$ such that $\chi(G)\le f(\omega(G))$ for every $H$-free graph $G$.
\end{thm}
This is open in general, but has been proved~\cite{gyarfas} when $H$ is a path, and for several other simple types of tree 
(\cite{distantstars, gyarfasprob, gst, kierstead, kierstead2, scott, newbrooms}; see~\cite{chibounded} for a survey).  The result is 
also known if all induced subdivisions of a tree are excluded~\cite{scott}.

A class of graphs is {\em hereditary} if the class is closed under taking induced subgraphs and under isomorphism, and a hereditary
class is said to be {\em $\chi$-bounded} if 
there is a function $f$ such that $\chi(G)\le f(\omega(G))$ for every graph $G$ in the class (thus the
Gy\'arf\'as-Sumner conjecture says that, for every forest $H$, the class of $H$-free graphs is $\chi$-bounded).  
Louis Esperet~\cite{esperet} made the following
conjecture:
\begin{thm}\label{Econj0}
	{\bf (False) Conjecture: }Let $\mathcal G$ be a $\chi$-bounded class.  Then there is a polynomial function $f$ such that
	$\chi(G)\le f(\omega(G))$ for every $G\in\mathcal G$.
\end{thm}
Esperet's conjecture was recently shown to be false by Bria\'nski, Davies and Walczak~\cite{brianski}.  However, this raises the further question: which $\chi$-bounded classes are polynomially $\chi$-bounded?
In particular, 
the two conjectures
\ref{GYconj} and \ref{Econj0} would together imply the following, which is still open:
\begin{thm}\label{Econj}
	{\bf Conjecture: }For every forest $H$, there exists $c>0$ such that $\chi(G)\le \omega(G)^c$ for every $H$-free graph $G$.
\end{thm}
This is a beautiful conjecture.  In most cases where the Gy\'arf\'as-Sumner conjecture has been proved, the current bounds are very far from polynomial, and \ref{Econj} has been only been proved for a much smaller collection of forests (see~\cite{liu, Kst, polystar,poly3, poly5,poly6,Schiermeyer}).
In~\cite{poly3} we proved it for any $P_5$-free tree $H$,
but it has not been settled for any tree $H$ that contains $P_5$.  In this paper we focus on the case $H=P_5$.

The best previously-known bound on the chromatic number of $P_5$-free graphs in terms of their clique number, 
due to Esperet, Lemoine, Maffray, and Morel~\cite{esperet2}, was exponential:
\begin{thm}\label{P5bound}
        If $G$ is $P_5$-free and $\omega(G)\ge 3$ then $\chi(G)\le (5/27)3^{\omega(G)}$.
\end{thm}
Here we make a significant improvement, showing  a ``near-polynomial'' bound:
\begin{thm}\label{mainthm}
        If $G$ is $P_5$-free and $\omega(G)\ge 3$ then $\chi(G)\le \omega(G)^{\log_2(\omega(G))}$.
\end{thm}
(The cycle of length five shows that we need to assume $\omega(G)\ge 3$. Sumner~\cite{sumner} showed that 
$\chi(G)\le 3$ when $\omega(G)=2$.)
Conjecture \ref{Econj} when $H=P_5$ is of great interest, because of a famous conjecture due to Erd\H{o}s and 
Hajnal~\cite{EH0, EH}, that:
\begin{thm}\label{EHconj}
        {\bf Conjecture: }For every graph $H$ there exists $c>0$ such that $\alpha(G)\omega(G)\ge |G|^c$ for every $H$-free graph $G$.
\end{thm}
This is open in general, despite a great deal of effort; and 
in view of~\cite{fivehole}, the smallest graph $H$ for which \ref{EHconj} is undecided is the graph $P_5$.
Every forest $H$ satisfying \ref{Econj} also
satisfies the Erd\H{o}s-Hajnal conjecture, and so showing that $H=P_5$ satisfies \ref{Econj} would be a significant result.

We use standard notation throughout.
When $X\subseteq V(G)$, $G[X]$ denotes the subgraph induced on $X$. We write $\chi(X)$ for $\chi(G[X])$
when there is no ambiguity.

\section{The main proof}
We denote the set of nonnegative real numbers by $\mathbb{R}_+$, and the set of nonnegative integers by $\mathbb{Z}_+$.
Let $f:\mathbb{Z}_+\rightarrow \mathbb{R}_+$ be a function. We say
\begin{itemize}
	\item $f$ is {\em non-decreasing} if
$f(y)\ge f(x)$ for all integers $x,y\ge 0$ with $y>x\ge 0$; 
	\item $f$ is a {\em binding} function for a graph $G$ if
it is non-decreasing and $\chi(H)\le f(\omega(H))$ for every induced subgraph $H$ of $G$; and 
	\item $f$ is a {\em near-binding} function for $G$ if $f$ is non-decreasing and 
$\chi(H)\le f(\omega(H))$ for every induced subgraph $H$ of $G$ different from $G$. 
\end{itemize}

In this section we show that if a function $f$ satisfies a certain inequality, then it is a binding function for all $P_5$-free 
graphs. Then at the end
we will give a function that satisfies the inequality, and deduce \ref{mainthm}.

A {\em cutset} in a graph $G$ is a set $X$ such that $G\setminus X$ is disconnected.
A vertex $v\in V(G)$ is {\em mixed} on a set $A\subseteq V(G)$ or a subgraph $A$ of a graph $G$ if $v$ is not in $A$ and has a neighbour and a non-neighbour in $A$. It is {\em complete}
to $A$ if it is adjacent to every vertex of $A$.
We begin with the following:
\begin{thm}\label{cutset}
	Let $G$ be $P_5$-free, and let $f$ be a near-binding function for $G$. Let $G$ be connected, and let $X$
	be a cutset of $G$. Then 
	$$\chi(G\setminus X)\le f(\omega(G)-1)+ \omega(G) f(\lfloor \omega(G)/2\rfloor).$$
\end{thm}
\Proof
We may assume (by replacing $X$ by a subset if necessary) that $X$ is a minimal cutset of $G$; and so $G\setminus X$ has at least two components, 
and every vertex in $X$ has a neighbour in $V(B)$, for every component $B$ of $G\setminus X$. Let $B$ be one such component; we will prove that
$\chi(B)\le f(\omega(G)-1)+ \omega(G) f(\lfloor \omega(G)/2\rfloor)$, from which the result follows.

Choose $v\in X$ (this is possible since $G$ is connected), and let $N$ be the set of vertices in $B$ adjacent to $v$. Let the components of $B\setminus N$
be $R_1\ll R_k, S_1\ll S_\ell$, where $R_1\ll R_k$ each have chromatic number more than $f(\lfloor \omega(G)/2\rfloor)$, 
and $S_1\ll S_\ell$ each have chromatic number
at most $f(\lfloor \omega(G)/2\rfloor)$. Let $S$ be the union of the graphs $S_1\ll S_\ell$; thus $\chi(S)\le f(\lfloor \omega(G)/2\rfloor)$. For $1\le i\le k$, let $Y_i$
be the set of vertices in $N$ with a neighbour in $V(R_i)$, and let $Y=Y_1\cupcup Y_k$.
\\
\\
(1) {\em For $1\le i\le k$, every vertex in $Y_i$ is complete to $R_i$.}
\\
\\
Let $y\in Y_i$. Thus $y$ has a neighbour in $V(R_i)$; suppose that $y$ is mixed on $R_i$. Since $R_i$ is connected, 
there is an edge $ab$
of $R_i$ such that $y$ is adjacent to $a$ and not to $b$. Now $v$ has a neighbour in each component of $G\setminus X$, and since there are 
at least two such components, there is a vertex $u\in V(G)\setminus (X\cup V(B))$ adjacent to $v$. But then $u\dd v\dd y\dd a\dd b$
is an induced copy of $P_5$, a contradiction. This proves (1).
\\
\\
(2) {\em $\chi(Y)\le (\omega(G)-1) f(\lfloor \omega(G)/2\rfloor)$.}
\\
\\
Let $1\le i\le k$. Since 
$f(\lfloor \omega(G)/2\rfloor)<\chi(R_i)\le f(\omega(R_i))$, and $f$ is non-decreasing, it follows that $\omega(R_i)>\omega(G)/2$. By (1), 
$\omega(G[Y_i])+\omega(R_i)\le \omega(G)$, and so 
$\omega(G[Y_i])<\omega(G)/2$. Consequently $\chi(Y_i)\le f(\lfloor \omega(G)/2\rfloor)$, for $1\le i\le k$. Choose $I\subseteq \{1\ll k\}$ minimal such that
$\bigcup_{i\in I}Y_i=Y$. From the minimality of $I$, for each $i\in I$ there exists $y_i\in Y_i$ such that
for each $j\in I\setminus \{i\}$ we have that $y_i\notin Y_j$; and so the vertices $y_i\;(i\in I)$ are all distinct. For each $i\in I$ choose $r_i\in V(R_i)$. 
For all distinct $i,j\in I$, if $y_i,y_j$ are nonadjacent, 
then $r_i\dd y_i\dd v\dd y_j\dd r_j$ is isomorphic to $P_5$, a contradiction. Hence the vertices $y_i\;(i \in I)$ are all pairwise adjacent, and adjacent to $v$;
and so
$|I|\le \omega(G)-1$. Thus $\chi(Y)=\chi(\bigcup_{i\in I}Y_i)\le (\omega(G)-1)f(\lfloor \omega(G)/2\rfloor)$. This proves (2).

\bigskip

All the vertices in $N\setminus Y$ are adjacent to $v$, and so $\omega(G[N\setminus Y])\le \omega(G)-1$. Moreover, for $1\le i\le k$, each vertex
of $R_i$ is adjacent to each vertex in $Y_i$, and $Y_i\ne \emptyset$ since $B$ is connected, and so $\omega(R_i)\le \omega(G)-1$.
Since there are no edges between any two of the graphs $G[N\setminus Y], R_1\ll R_k$, their union ($Z$ say) has clique number at most $\omega(G)-1$ and so 
has chromatic number at most $f(\omega(G)-1)$. But $V(B)$ is the union of $Y, V(S)$ and $V(Z)$; and so 
$$\chi(B)\le f(\omega(G)-1)+(\omega(G)-1)f(\lfloor \omega(G)/2\rfloor) +f(\lfloor \omega(G)/2\rfloor).$$
This proves \ref{cutset}.~\bbox

\begin{thm}\label{ineqs}
	Let $\Omega\ge 1$, and let $f:\mathbb{Z}_+\rightarrow \mathbb{R}_+$ be non-decreasing, satisfying the following:
	\begin{itemize}
		\item $f$ is a binding function for every $P_5$-free graph $H$ with $\omega(H)\le \Omega$; and
		\item $f(w-1) + (w+2) f(\lfloor w/2\rfloor) \le f(w)$ for each integer $w> \Omega$.
	\end{itemize}
	Then $f$ is a binding function for every $P_5$-free graph $G$.
\end{thm}
\Proof
We prove by induction on $|G|$ that if $G$ is $P_5$-free then $f$ is a binding function for $G$. Thus, we may assume that $G$ is $P_5$-free and $f$ is near-binding for $G$. If $G$
is not connected, or $\omega(G)\le \Omega$, it follows that $f$ is binding for $G$, so we assume that $G$ is connected and $\omega(G)>\Omega$.
Let us write $w=\omega(G)$ and $m=\lfloor w/2\rfloor$. If $\chi(G)\le f(w)$ then $f$ is a binding function for $G$, so we assume, for a contradiction, that:
\\
\\
(1) {\em $\chi(G)>f(w-1) + (w+2)f(m)$.}
\\
\\
We deduce that:
\\
\\
(2) {\em Every cutset $X$ of $G$ satisfies $\chi(X)> 2f(m)$.}
\\
\\
If some cutset $X$ satisfies $\chi(X)\le 2f(m)$, then since $\chi(G\setminus X)\le f(w-1)+ w f(m)$ by \ref{cutset},
it follows that 
$\chi(G)\le f(w-1)+(w+2)f(m)$, contrary to (1).
This proves (2).
\\
\\
(3) {\em If $P,Q$ are cliques of $G$, both of cardinality at least $w/2$, then $G[P\cup Q]$ is connected.}
\\
\\
Suppose not; then there is a minimal subset $X\subseteq V(G)\setminus (P\cup Q)$ such that $P,Q$ are subsets of different components
($A,B$ say) of $G\setminus X$. From the minimality of $X$, every vertex $x\in X$ has a neighbour in $V(A)$ and a neighbour in $V(B)$. 
If $x$ is mixed on $A$ and mixed on $B$, then since $A$ is connected, there is an edge $a_1a_2$ of $A$ such that $x$ is adjacent 
to $a_1$ and not to $a_2$;
and similarly there is an edge $b_1b_2$ of $B$ with $x$ adjacent to $b_1$ and not to $b_2$. But 
then $a_2\dd a_1\dd x\dd b_1\dd b_2$ is an induced copy of $P_5$, a contradiction; so every $x\in X$ is complete to at least one of $A,B$.
The set of vertices in $X$ complete to $A$ is also complete to $P$, and hence has clique number at most $m$, and hence has chromatic
number at most $f(m)$; and the same for $B$. Thus $\chi(X)\le 2f(m)$, contrary to (2). This proves (3).

\bigskip

If $v\in V(G)$, we denote its set of neighbours by $N(v)$, or $N_G(v)$.
Let $a\in V(G)$, and let $B$ be a component of $G\setminus (N(a)\cup \{a\})$;
we will show that $\chi(B)\le (w-m+2)f(m)$. 

A subset $Y$ of $V(B)$ is a {\em joint} of $B$ if there is a component $C$ of $B\setminus Y$ such that $\chi(C)> f(m)$
and $Y$ is complete to $C$. If $\emptyset$ is not a joint of $B$ then $\chi(B)<f(m)$ and the claim holds, so
we may assume that $\emptyset$ is a joint of $B$; let $Y$ be a joint of $B$ chosen with $Y$ maximal, and let 
$C$ be a component of $B\setminus Y$ such that $\chi(C)> f(m)$
and $Y$ is complete to $C$.
\\
\\
(4) {\em If $v\in N(a)$ has a neighbour in $V(C)$, then $\chi(V(C)\setminus N(v))\le f(m)$.}
\\
\\
Let $N_C(v)$ be the set of neighbours of $v$ in $V(C)$, and $M=V(C)\setminus N_C(v)$; and suppose that $\chi(M)> f(m)$. Let
$C'$ be a component of $G[M]$
with $\chi(C')> f(m)$, and let $Z$ be the set of vertices in $N_C(v)$ that have a neighbour in $V(C')$. Thus $Z\ne \emptyset$, 
since 
$N_C(v),V(C')\ne \emptyset$ and $C$ is connected. If some $z\in Z$ is mixed on $C'$, let $p_1p_2$ be an edge of $C'$ such that
$z$ is adjacent to $p_1$ and not to $p_2$; then $a\dd v\dd z\dd p_1\dd p_2$ is an induced copy of $P_5$, a contradiction. So 
every vertex in $Z$ is complete to $V(C')$; but also every vertex in $Y$ is complete to $V(C)$ and hence to $V(C')$, and so 
$Y\cup Z$ is a joint of $B$, contrary to the maximality of $Y$. This proves (4).
\\
\\
(5) {\em $\chi(Y)\le f(m)$ and $\chi(C)\le (w-m+1) f(m)$.}
\\
\\
Let $X$ be the set of vertices in $N(a)$ that have a neighbour in $V(C)$. Since $C$ is a component of $B\setminus Y$
and hence a component of $G\setminus (X\cup Y)$, and $a$ belongs to a different component of $G\setminus (X\cup Y)$, it follows
that $X\cup Y$ is a cutset of $G$.
By (2), $\chi(X\cup Y)>2f(m)$. Since $\omega(C)\ge m+1$ (because $\chi(C)>f(m)$, and $f$ is near-binding for $G$)
and every vertex in $Y$ is complete to $V(C)$, it follows that $\omega(G[Y])\le w-m-1\le m$, and so has chromatic number
at most $f(m)$ as claimed; and so $\chi(X)>f(m)$. Consequently there is a clique $P\subseteq X$ with 
cardinality $w-m$.
The subgraph induced on the set of vertices of $C$ complete to $P$ has clique number at most $m$, and so has chromatic number at most 
$f(m)$; and for each $v\in P$, the set of vertices of $C$ nonadjacent to $v$ has chromatic number at most 
$f(m)$ by (4). Thus 
$\chi(C)\le (|P|+1)f(m)= (w-m+1)f(m)$. This proves (5).
\\
\\
(6) {\em $\chi(B)\le (w-m+2) f(m)$.}
\\
\\
By (3), every clique contained in $V(B)\setminus (V(C)\cup Y)$ has cardinality less 
than $w/2$ (because it is anticomplete to the largest clique of $C$) and so 
$$\chi(B\setminus (V(C)\cup Y))\le f(m);$$ 
and hence $\chi(B\setminus Y)\le (w-m+1)f(m)$ by (5), since
there are no edges between $C$ and $V(B)\setminus (V(C)\cup Y)$. But $\chi(Y)\le f(m)$ by (5), and so 
$\chi(B)\le (w-m+2) f(m)$. This proves (6).

\bigskip

By (6), $G\setminus N(a)$ has chromatic number at most $(w-m+2)f(m)$. 
But $G[N(a)]$ has clique number at most $w-1$ and so chromatic number
at most $f(w-1)$; and so $\chi(G)\le f(w-1)+(w-m+2)f(m)$, contrary to (1).  This proves \ref{ineqs}.~\bbox

Now we deduce \ref{mainthm}, which we restate:
\begin{thm}\label{mainthm2}
If $G$ is $P_5$-free and $\omega(G)\ge 3$ then $\chi(G)\le \omega(G)^{\log_2(\omega(G))}$.
\end{thm}
\Proof
Define $f(0)=0$, $f(1)=1$, $f(2)=3$, and $f(x)=x^{\log_2(x)}$ for every real number $x\ge 3$. Let $G$ be $P_5$-free. If 
$\omega(G)\le 2$ then $\chi(G)\le 3=f(2)$, by a result of Sumner~\cite{sumner}; if $\omega(G) = 3$ then $\chi(G)\le 5\le f(3)$,
 by an application of the result \ref{P5bound} of Esperet, Lemoine, Maffray, and Morel~\cite{esperet2};
and if $\omega(G)=4$ then $\chi(G)\le 15\le f(4)$, by another application of \ref{P5bound}. 
Consequently  
every $P_5$-free graph $G$ with clique number at most four has chromatic 
number at most $f(\omega(G))$. 

We claim that 
$$f(x-1) + (x+2)f(\lfloor x/2 \rfloor) \le f(x)$$ 
for each integer $x>4$.  If that is true, then 
by \ref{ineqs} with $\Omega=4$, we deduce that $\chi(G)\le f(\omega(G))$ for every $P_5$-free
graph $G$, and so \ref{mainthm} holds. Thus, it remains to show that 
$$f(x-1) + (x+2)f(\lfloor x/2 \rfloor) \le f(x)$$
for each integer $x>4$.  This can be verified by direct calculation when $x=5$, so we may assume that $x\ge 6$. 

The derivative of $f(x)/x^4$ is 
$$(2\log_2(x)-4)x^{\log_2(x)-5},$$
and so is  nonnegative for $x\ge 4$. Consequently 
$$\frac{f(x-1)}{(x-1)^4}\le \frac{f(x)}{x^4}$$
for $x\ge 5$. Since $x^2(x^2-2x-4)\ge (x-1)^4$ when $x\ge 5$, it follows that
$$\frac{f(x-1)}{x^2-2x-4}\le \frac{f(x)}{x^2},$$
that is, 
$$f(x-1) + \frac{2x+4}{x^2}f(x) \le f(x),$$
when $x\ge 5$.
But when $x\ge 6$ (so that $f(x/2)$ is defined and the first equality below holds), we have
$$f(\lfloor x/2\rfloor)\le f(x/2)=(x/2)^{\log_2(x/2)} = (x/2)^{\log_2(x) - 1}=(2/x)(x/2)^{\log_2(x)} = (2/x^2)f(x),$$
and so
$$f(x-1) + (x+2)f(\lfloor x/2\rfloor) \le f(x)$$
when $x\ge 6$. This proves \ref{mainthm2}.~\bbox

\end{document}